\documentclass{amsart}
 \newtheorem{thm}{Theorem}[section]
 \newtheorem{cor}[thm]{Corollary}
 \newtheorem{lem}[thm]{Lemma}
 \newtheorem{prop}[thm]{Proposition}
 \theoremstyle{definition}
 \newtheorem{defn-new}[thm]{Definition}
 \theoremstyle{remark}
 \newtheorem{rem-new}[thm]{Remark}
 \newtheorem*{ex-new}{Example}
 \numberwithin{equation}{section}

\begin{document}

%

\title{Generalized $(\kappa,\mu)$-space forms}

\author[A. Carriazo]{Alfonso Carriazo}
\address{Department of Geometry and Topology, Faculty of Mathematics, University of Seville,
           Apdo. de Correos 1160, 41080  -- Sevilla, SPAIN.}
\email{carriazo@us.es}

\author[V. Mart\'{\i}n Molina]{Ver\'onica Mart\'{\i}n Molina}
\address{Department of Geometry and Topology, Faculty of Mathematics, University of Seville,
           Apdo. de Correos 1160, 41080  -- Sevilla, SPAIN.}
\email{veronicamartin@us.es}

\thanks{The first two authors are partially supported by the MICINN grant MTM2011-22621 and the  PAI group FQM-327 (Junta de Andaluc\'{\i}a, Spain). The second one is also supported by the FPU program of the Ministerio de Educaci\'on, Spain.}
\author[M. M. Tripathi]{Mukut Mani Tripathi}
\address{Department of Mathematics, Banaras Hindu University, Varanasi, 221 005, INDIA}
\email{mmtripathi66@yahoo.com}
%


\date{}

\begin{abstract}
Generalized $\left( \kappa ,\mu \right) $-space forms are
introduced and studied. We examine in depth the contact metric case and
present examples for all possible dimensions. We also analyse the
trans-Sasakian case.
\end{abstract}

\maketitle

\section{Introduction}\label{sect-intro}

A generalized Sasakian space form was defined by the first named
author (jointly with P. Alegre and D. E. Blair) in \cite{ABC}
as that almost contact metric manifold $( M^{2n+1},\phi ,\xi,\eta ,g) $ whose curvature tensor $R$ is given by
\begin{equation}
R=f_{1}R_{1}+f_{2}R_{2}+f_{3}R_{3},  \label{eq-R-gssf}
\end{equation}
where $f_{1},f_{2},f_{3}$ are some differentiable functions on $M$
and
\begin{align*}
R_{1}( X,Y) Z &=g( Y,Z) X-g( X,Z) Y, \\
R_{2}( X,Y) Z &=g( X,\phi Z) \phi Y-g(Y,\phi Z)\phi X+2g(X,\phi Y)\phi Z, \\
R_{3}( X,Y) Z &=\eta ( X) \eta (Z) Y-\eta ( Y) \eta ( Z) X+g(X,Z) \eta (Y) \xi -g( Y,Z) \eta ( X) \xi ,
\end{align*}
for any vector fields $X,Y,Z$ on $M$. We denote it by $M(f_{1},f_{2},f_{3})$.

Since then, several papers have appeared concerning different
aspects of this topic: structures, submanifolds and conformal
changes of metric (\cite{AC-DGA}, \cite{AC-Taiwanese} and
\cite{AC-preprint}), B.-Y. Chen's inequalities (\cite{ACKY}),
slant submanifolds inheriting the structure (\cite{ACOS}),
CR-submanifolds (\cite{Gh-So-Sh-06-JP} and
\cite{Gh-So-Sh-06-BJGA}), the Ricci curvature of some submanifolds
(\cite{hong}, \cite{mihai}  and \cite{st}), conformal flatness and
local symmetry (\cite{KimUK-06}), some other symmetry properties
(\cite{ggb} and \cite{gkb}) and immersions of warped products
(\cite{olt} and \cite{Yoon-04-IMJ}). Other related papers are \cite{AyL1} and \cite{AyL2}.

On the other hand, a contact metric manifold $(M^{2n+1},\phi ,\xi ,\eta ,g)$ is said to be a generalized $\left( \kappa ,\mu \right) $-space
if its curvature tensor satisfies the condition
\[
R\left( X,Y\right) \xi =\kappa \left\{ \eta \left( Y\right) X-\eta
\left( X\right) Y\right\} +\mu \left\{ \eta \left( Y\right) hX-\eta
\left( X\right) hY\right\},
\]
for some smooth functions $\kappa $ and $\mu $ on $M$ independent of
the choice of vectors fields $X$ and $Y$. If $\kappa $ and $\mu $
are constant, the manifold is called a $\left( \kappa ,\mu \right)
$-space. T.
Koufogiorgos proved in \cite{Kouf-97} that if a $\left( \kappa ,\mu \right)$-space $M$ has constant $\phi $-sectional curvature $c$ and dimension greater than 3, the curvature tensor of this $\left( \kappa ,\mu \right)$-space form is given by
\begin{equation}
R=\frac{c+3}{4}\,R_{1}+\frac{c-1}{4}\,R_{2}+\left(
\frac{c+3}{4}-\kappa \right) R_{3}+R_{4}+\frac{1}{2}R_{5}+\left(
1-\mu \right) \,R_{6}, \label{eq-R-(k,m)sf}
\end{equation}
where $R_{1},R_{2},R_{3}$ are the tensors defined above and
\begin{align*}
R_{4}(X,Y)Z=&g( Y,Z) hX-g( X,Z)hY+g( hY,Z) X-g(hX,Z)Y, \\
R_{5}(X,Y)Z =&g( hY,Z) hX-g(hX,Z)hY\\
& +g(\phi hX,Z)\phi hY-g(\phi hY,Z)\phi hX, \\
R_{6}(X,Y)Z=&\eta(X)\eta(Z)hY-\eta(Y)\eta(Z)hX\\
             &  +g(hX,Z)\eta(Y)\xi-g(hY,Z)\eta(X)\xi,
\end{align*}
for any vector fields $X,Y,Z$, where $2h=L_{\xi }\phi $ and $L$\ is
the usual Lie derivative. Now, a natural question arises: is it
possible to generalize the notion of $\left( \kappa ,\mu \right)
$-space form, by replacing the constants in equation
(\ref{eq-R-(k,m)sf}) by some differentiable functions?

Therefore, in view of the previous works on generalized Sasakian
space forms, we say that an almost contact metric manifold
$(M^{2n+1},\phi ,\xi ,\eta ,g)$ is a {\em generalized} $\left(
\kappa ,\mu \right) ${\em -space
form} if there exist differentiable functions $f_{1},f_{2},f_{3}$, $f_{4},f_{5},f_{6}$ on $M$ such that
\[
R=f_{1}R_{1}+f_{2}R_{2}+f_{3}R_{3}+f_{4}R_{4}+f_{5}R_{5}+f_{6}R_{6}.
\]It is obvious that $\left( \kappa ,\mu \right) $-space forms are natural
examples of generalized $\left( \kappa ,\mu \right) $-space forms,
with constant functions
\[
f_{1}=\frac{c+3}{4},\;f_{2}=\frac{c-1}{4},\;f_{3}=\frac{c+3}{4}-\kappa
,\;f_{4}=1,\;f_{5}=\frac{1}{2},\;f_{6}=1-\mu .
\]
We also have generalized Sasakian space forms, with
$f_{4}=f_{5}=f_{6}=0$.

Thus, in this paper we introduce and study generalized $\left(
\kappa ,\mu
\right) $-space forms. The paper is organized as follows. The section~\ref{sect-prelim} contains some necessary background on almost contact
metric geometry. Afterwards, in section~\ref{sect-def} we formally
give the definition of generalized $\left( \kappa ,\mu \right)
$-space form and check that some results that were true for
generalized Sasakian space forms are also correct for generalized
$\left( \kappa ,\mu \right) $-space forms. Then, we obtain some
basic identities for generalized $\left( \kappa ,\mu \right)$-space
forms which are analogous to those satisfied by Sasakian manifolds.
In section~\ref{sect-contact} we prove that contact metric
generalized $\left( \kappa ,\mu \right) $-space forms are generalized $\left( \kappa ,\mu \right) $-spaces with $\kappa =f_{1}-f_{3}$ and
$\mu =f_{4}-f_{6}$. Next, we observe that if dimension is greater
than or equal
to $5$, then they are $(-f_{6},1-f_{6})$-spaces with constant $\phi$-sectional curvature $2f_{6}-1$. Furthermore, $f_{4}=1$, $f_{5}=1/2$ and $f_{1},f_{2},f_{3}$ depend linearly on the constant $f_{6}$. We also
give a method for constructing infinitely many examples of this
type. Later, we pay attention to the $3$-dimensional case, in
which we prove that the expression for the curvature tensor is not
unique and that several properties and results must be satisfied. We
also check that the example of generalized $\left( \kappa ,\mu
\right) $-space with non-constant $\kappa $ and $\mu $ that T.
Koufogiorgos and C. Tsichlias provided in \cite{Kouf-00} is a
generalized $\left( \kappa ,\mu \right) $-space form with
non-constant
functions $f_{1},f_{3}$ and $f_{4}$. Finally, we prove in section~\ref{sect-trans} that if a manifold is trans-Sasakian, then $h=0$.
Therefore, generalized $\left( \kappa ,\mu \right) $-space forms
with trans-Sasakian
structure are generalized Sasakian space forms, already studied in \cite{AC-DGA}.

\section{Preliminaries\label{sect-prelim}}

In this section, we recall some general definitions and basic
formulas which will be used later. For more background on almost
contact metric manifolds, we recommend the reference
\cite{Blair-02}. Anyway, we will recall some more specific notions
and results in the following sections, when required.

An odd-dimensional Riemannian manifold $(M,g)$ is said to be an {\em
almost
contact metric manifold} if there exist on $M$ a $(1,1)$-tensor field $\phi $, a vector field $\xi $ (called the {\em structure vector
field}) and a 1-form $\eta $ such that $\eta (\xi )=1$, $\phi
^{2}X=-X+\eta \left( X\right) \xi $ and $g(\phi X,\phi Y)=g\left(
X,Y\right) -\eta \left( X\right) \eta \left( Y\right) $ for any
vector fields $X,Y$ on $M$. In particular, on an almost contact
metric manifold we also have $\phi \xi =0$ and $\eta \circ \phi =0$.

Such a manifold is said to be a {\em contact metric manifold} if ${\rm d} \eta =\Phi $, where $\Phi \left( X,Y\right) =g(X,\phi Y)$ is the {\em fundamental} $2${\em -form} of $M$. If, in addition, $\xi $ is a
Killing vector field, then $M$ is said to be a $K$-{\em contact
manifold}. It is well-known that a contact metric manifold is a
$K$-contact manifold if and only if
\begin{equation}
\nabla _{X}\xi =-\,\phi X  \label{eq-K-cont}
\end{equation}
for all vector fields $X$ on $M$. Even an almost contact metric
manifold satisfying the equation (\ref{eq-K-cont}) becomes a
$K$-contact manifold.

On the other hand, the almost contact metric structure of $M$ is
said to be {\em normal} if the Nijenhuis torsion $[\phi ,\phi ]$\ of
$\phi $ equals $-2{\rm d}\eta \otimes \xi $. A normal contact metric
manifold is called a {\em Sasakian manifold}. It can be proved that
an almost contact metric manifold is Sasakian if and only if
\begin{equation}
(\nabla _{X}\phi )Y=g\left( X,Y\right) \xi -\eta \left( Y\right) X
\label{eq-Sas}
\end{equation}
for any vector fields $X,Y$ on $M$. Moreover, for a Sasakian
manifold the following equation holds:
\[
R\left( X,Y\right) \xi =\eta \left( Y\right) X-\eta \left( X\right)
Y.
\]

Given an almost contact metric manifold $(M,\phi ,\xi ,\eta ,g)$, a $\phi $-{\em section} of $M$ at $p\in M$ is a section $\Pi \subseteq T_{p}M $ spanned by a unit vector $X_{p}$ orthogonal to $\xi _{p}$, and $\phi X_{p}$. The $\phi $-{\em sectional curvature of} $\Pi $ is
defined by $K(X,\phi X)=R(X,\phi X,\phi X,X)$. A Sasakian manifold
with constant $\phi $-sectional curvature $c$ is called a {\em
Sasakian space form}. In such a case, its Riemann curvature tensor
is given by equation (\ref{eq-R-gssf}) with functions $f_{1} =
(c+3)/4$ and $f_{2} = f_{3}= (c-1)/4$.

It is well known that on a contact metric manifold $\left( M,\phi ,\xi ,\eta ,g\right) $, the tensor $h$, defined by $2h=L_\xi \phi$, is symmetric and satisfies the following relations \cite{Blair-02}
\begin{equation}
h\xi =0,\quad \nabla _{X}\xi =-\phi X-\phi hX,\quad h\phi =-\phi
h,\quad {\rm tr}(h)=0,\quad \eta \circ h=0. \label{eq-cont-h}
\end{equation}
Therefore, it follows from equations (\ref{eq-K-cont}) and
(\ref{eq-cont-h}) that a contact metric manifold is $K$-contact if
and only if $h=0$.

Finally, we assume that all the functions considered in this paper
will be differentiable functions on the corresponding manifolds.

\section{Definition and first results\label{sect-def}}

In this section we give the formal definition of generalized $\left(
\kappa ,\mu \right) $-space forms and prove some basic results about
these manifolds. Then we study some interesting properties of their
curvature tensor.

\begin{defn-new}
\label{defn-gkmsf} We say that an almost contact metric manifold
$(M,\phi ,\xi ,\eta ,g)$ is a {\em generalized $\left( \kappa ,\mu
\right) $-space form} if there exist functions
$f_{1},f_{2},f_{3},f_{4},f_{5},f_{6}$ defined on $M$ such that
\begin{equation}
R=f_{1}R_{1}+f_{2}R_{2}+f_{3}R_{3}+f_{4}R_{4}+f_{5}R_{5}+f_{6}R_{6},
\label{eq-R-gkmsf}
\end{equation}
where $R_{1},R_{2},R_{3},R_{4},R_{5},R_{6}$ are the following
tensors
\begin{align*}
R_{1}( X,Y) Z =&g( Y,Z) X-g( X,Z) Y, \\
R_{2}( X,Y) Z =&g( X,\phi Z) \phi Y-g(Y,\phi Z)\phi X+2g(X,\phi Y)\phi Z, \\
R_{3}( X,Y) Z =&\eta ( X) \eta (Z) Y-\eta ( Y) \eta ( Z) X+g(X,Z) \eta (Y) \xi -g( Y,Z) \eta ( X) \xi , \\
R_{4}( X,Y) Z =&g( Y,Z) hX-g( X,Z)hY+g( hY,Z) X-g(hX,Z)Y, \\
R_{5}( X,Y) Z =&g( hY,Z) hX-g(hX,Z)hY\\
&+g(\phi hX,Z)\phi hY-g(\phi hY,Z)\phi hX, \\
R_{6}( X,Y) Z =&\eta ( X) \eta (Z) hY-\eta ( Y) \eta ( Z) hX\\
&+g(hX,Z)\eta( Y) \xi -g( hY,Z) \eta ( X) \xi,
\end{align*}
for all vector fields $X,Y,Z$ on $M$, where $2h=L_{\xi }\phi $ and
$L$ is the usual Lie derivative. We will denote such a manifold by
$M(f_{1},\ldots ,f_{6})$.
\end{defn-new}

\begin{rem-new}
It is obvious that $\left( \kappa ,\mu \right) $-contact space forms
of dimension greater than $3$ are natural examples of generalized
$\left( \kappa ,\mu \right) $-space forms, where
\[
f_{1}=\frac{c+3}{4},\;f_{2}=\frac{c-1}{4},\;f_{3}=\frac{c+3}{4}-\kappa
,\;f_{4}=1,\;f_{5}=\frac{1}{2},\;f_{6}=1-\mu
\]
are constant. Generalized Sasakian space forms (defined in
\cite{ABC}) are also examples with $f_{4}=f_{5}=f_{6}=0$ and
$f_{1},f_{2},f_{3}$ not necessarily constant.
\end{rem-new}

As we have already pointed out, $h=0$ for a $K$-contact manifold.
Therefore, a generalized $\left( \kappa ,\mu \right) $-space form
with such a structure is actually a generalized Sasakian space form.
Hence, the following results
are inferred from Proposition 3.6, Theorem~3.7 and Theorem~3.15 from \cite{ABC}:

\begin{thm}
\label{th-Sas-1} Let $M(f_{1},\ldots ,f_{6})$ be a generalized
$\left(
\kappa ,\mu \right) $-space form. If $M$ is a $K$-contact manifold, then $f_{3}=f_{1}-1$. Moreover, $M$ is Sasakian.
\end{thm}

\begin{thm}
\label{th-Sas-2} Let $M(f_{1},\ldots ,f_{6})$ be a generalized
$\left(
\kappa ,\mu \right) $-space form. If $M$ is a Sasakian manifold, then $f_{2}=f_{3}=f_{1}-1$.
\end{thm}

\begin{rem-new}
Sasakian manifolds are always $K$-contact, while the converse is not
true in general, only in dimension 3. However, we have just seen
that being Sasakian and $K$-contact are equivalent concepts for
generalized $\left( \kappa ,\mu \right) $-space forms.
\end{rem-new}

Using the properties of $h$, it can be proved that:

\begin{thm}
\label{th-Sas-3} Let $M(f_{1},\ldots ,f_{6})$ be a generalized
$\left(
\kappa ,\mu \right) $-space form. If $M$ is a contact metric manifold with $f_{3}=f_{1}-1$, then it is a Sasakian manifold.
\end{thm}

\begin{proof} Let $M^{2n+1}$ be a contact metric manifold
satisfying $f_3=f_1-1$. Because of Theorem~\ref{th-Sas-1} we would
only need
to prove that $M$ is $K$-contact, which is equivalent to checking that $S(\xi,\xi)=2n$, where $S$ denotes the Ricci curvature tensor (see \cite{Blair-02}, p.~92).

If we take a local orthonormal basis $\{e_{1},\ldots ,e_{2n},\xi
\}$, then a direct computation from (\ref{eq-R-gkmsf}) gives
\[
R(e_{i},\xi ,\xi ,e_{i})=1+(f_{4}-f_{6})g(he_{i},e_{i}),
\]
where we have used the properties of almost contact metric
manifolds, the
fact that $h\xi =0$ and the hypothesis $f_{1}-f_{3}=1$. Therefore,
\begin{align*}
S\left( \xi ,\xi \right) & = {\sum_{i=1}^{2n}R}\left( {e_{i},\xi ,\xi ,e_{i}}\right) +R(\xi ,\xi ,\xi ,\xi ) =
\sum_{i=1}^{2n}(1+(f_{4}-f_{6})g(he_{i},e_{i})) \\
& =  2n + (f_{4} - f_{6})\sum_{i=1}^{2n}g(he_{i},e_{i}) = 2n+(f_{4}-f_{6}) {\rm tr}(h) = 2n,
\end{align*}
because of (\ref{eq-cont-h}).
\end{proof}

We will now calculate $K(X,\phi X)$ (the $\phi $-sectional
curvature), $K(X,\xi )$ (the $\xi $-sectional curvature) and $K(\phi
X,\xi )$ for generalized $\left( \kappa ,\mu \right) $-space forms.

\begin{prop}\label{prop-gkmsf-K(X,phiX)}
Let $M(f_{1},\ldots ,f_{6})$ be a generalized $\left( \kappa ,\mu \right) $-space form. Then the $\phi $-sectional
curvature of the $\phi $-section spanned by the unit vector field
$X$, orthogonal to $\xi $, is given by
\[
K(X,\phi X)=f_{1}+3f_{2}+f_{4}\,g((h-\phi h\phi )X,X).
\]
If $M$ is also a contact metric manifold, then the $\phi $-sectional
curvature is given by
\[
K(X,\phi X)=f_{1}+3f_{2},
\]
so it does not depend on the choice of vector field $X$.
\end{prop}

\begin{proof} By applying the properties of an almost
contact
metric manifold, we can directly calculate from (\ref{eq-R-gkmsf}) the $\phi $-sectional curvature of the $\phi $-section spanned by $\{X,\phi X\}$, where $X$ is a unit vector field orthogonal to $\xi$, as
follows:
\begin{align*}
K\left( X,\phi X\right) &=f_{1}+3f_{2}+f_{4}\left\{ g\left(hX,X\right)
+g\left( h\phi X,\phi X\right) \right\} \\
&=f_{1}+3f_{2}+f_{4}g\left( \left( h-\phi h\phi \right) X,X\right).
\end{align*}
If $M$ is a contact metric manifold then, in view of (\ref{eq-cont-h}), $h\phi =-\phi h$ and we get
\[
\left( h-\phi h\phi \right) X=\left( h+h\phi ^{2}\right) X=h\left(
X+\phi ^{2}X\right) =h\left( X-X\right) =0.
\]
Therefore, in such a case, $K\left( X,\phi X\right) =f_{1}+3f_{2}$.
\end{proof}

A direct computation, similar to that of Theorem~\ref{th-Sas-3},
gives:

\begin{prop}
\label{prop-gkmsf-K(X,xi)} Let $M(f_{1},\ldots ,f_{6})$ be a generalized $\left( \kappa ,\mu \right) $-space form. Then the $\xi $-sectional
curvature
of the $\xi$-section spanned by the unit vector field $X$, orthogonal to $\xi $, is given by
\[
K(X,\xi )=f_{1}-f_{3}+(f_{4}-f_{6})g(hX,X).
\]
\end{prop}

\begin{cor}
\label{cor-gkmsf-K(phiX,xi)} Let $M(f_{1},\ldots ,f_{6})$ be a generalized $\left( \kappa ,\mu \right) $-space form. If $X$ is a unit vector
field orthogonal to $\xi$, then
\[
K(\phi X,\xi )=f_{1}-f_{3}+(f_{4}-f_{6})g(h\phi X,\phi X).
\]
If $M$ is also a contact metric manifold, then
\[
K(\phi X,\xi )=f_{1}-f_{3}-(f_{4}-f_{6})g(hX,X).
\]
\end{cor}

\begin{proof} The first assertion follows directly from Proposition \ref{prop-gkmsf-K(X,xi)}. With respect to the contact metric case,
we just have to take into account that
\[
g(h\phi X,\phi X)=-g(\phi hX,\phi X)=-g(hX,X),
\]
by virtue of (\ref{eq-cont-h}).
\end{proof}

The following classic result appears in (\cite{Blair-02}, pp.
94-95):

\begin{lem}
\label{lem-Sas-1} Let $M$ be a Sasakian manifold. If we put
\begin{align*}
\widetilde{P}\left( X,Y,Z,W\right) &=d\eta \left( X,Z\right) g(Y,W)-d\eta (X,W)g\left( Y,Z\right) \\
&-\ d\eta \left( Y,Z\right) g(X,W)+d\eta (Y,W)g(X,Z),
\end{align*}
then
\begin{equation}
R(X,Y,Z,\phi W)+R(X,Y,\phi Z,W)=-\ \widetilde{P}\left(
X,Y,Z,W\right)
\end{equation}
for any vector fields $X,Y,Z,W$ on $M$ and
\begin{equation}
R(\phi X,\phi Y,\phi Z,\phi W)=R\left( X,Y,Z,W\right) ,
\end{equation}
\begin{equation}
R(X,\phi X,Y,\phi Y)=R(X,Y,X,Y)+R(X,\phi Y,X,\phi Y)-2\widetilde{P}(X,Y,X,\phi Y)  \label{eq-Blair-1}
\end{equation}
for any vector fields $X,Y,Z,W$ orthogonal to $\xi $.
\end{lem}

Let $\left( M,\phi ,\xi ,\eta ,g\right) $ be any almost contact
metric manifold. We now denote
\begin{align*}
P\left( X,Y,Z,W\right) &=g\left( X,\phi Z\right) g(Y,W)-g(X,\phi
W)g\left( Y,Z\right) \\
&-\ g(Y,\phi Z)g(X,W)+g(Y,\phi W)g\left( X,Z\right)
\end{align*}
for any vectors fields $X,Y,Z,W$ on $M$. In particular, if $M$ is a
contact metric manifold, $P=\widetilde{P}$. We will study whether
similar results hold true for generalized $\left( \kappa ,\mu
\right) $-space forms.

We omit the corresponding proofs because they can be easily obtained
by
making direct computations from (\ref{eq-R-gkmsf}) and using (\ref{eq-cont-h}) in the contact metric cases.

\begin{prop}
Let $M(f_{1},\ldots ,f_{6})$ be a generalized $\left( \kappa ,\mu \right)$-space form. Given $X,Y,Z,W$ orthogonal to $\xi $, we have
\[
R_{i}(\phi X,\phi Y,\phi Z,\phi W)=R_{i}\left( X,Y,Z,W\right) \quad
{\rm for\quad }i=1,2,3,6.
\]
Therefore,
\begin{align*}
R(\phi X,\phi Y,\phi Z,\phi W)-f_{4}R_{4}(\phi X,\phi Y,\phi Z,\phi W)-f_{5}R_{5}(\phi X,\phi Y,\phi Z,\phi W) \\
\qquad =\ R\left( X,Y,Z,W\right) -f_{4}R_{4}\left( X,Y,Z,W\right) -f_{5}R_{5}\left( X,Y,Z,W\right).
\end{align*}
If $M$ is also a contact metric manifold, then
\[
R_{4}(\phi X,\phi Y,\phi Z,\phi W)=-\ R_{4}\left( X,Y,Z,W\right),
\]
\[
R_{5}(\phi X,\phi Y,\phi Z,\phi W)=R_{5}\left( X,Y,Z,W\right),
\]
and therefore
\[
R(\phi X,\phi Y,\phi Z,\phi W)=R\left( X,Y,Z,W\right)
-2f_{4}R_{4}\left( X,Y,Z,W\right) .
\]
\end{prop}

\begin{prop}
Let $M(f_{1},\ldots ,f_{6})$ be a generalized $\left( \kappa ,\mu \right) $-space form. Then
\begin{align*}
R(X,Y,& Z,\phi W)+R(X,Y,\phi Z,W)=-(f_{1}-f_{2})P(X,Y,Z,W)\\
&-f_{4}(P(hX,Y,Z,W)+P(X,hY,Z,W))-2f_{5}P(hX,hY,Z,W)
\end{align*}
for any vector fields $X,Y,Z,W$ orthogonal to $\xi $.
\end{prop}

It is clear that the above results extend Propositions 3.14 and 3.17
from \cite{ABC}. To obtain a result similar to equation
(\ref{eq-Blair-1}), we need the manifold to be a contact metric one:

\begin{prop}
Let $M(f_{1},\ldots ,f_{6})$ be a generalized $\left( \kappa ,\mu \right)$-space form. If $M$ is also a contact metric manifold, then
\begin{align*}
R(X,\phi X,Y,\phi Y) =& R(X,Y,X,Y)+R(X,\phi Y,X,\phi Y) \\
& -\ 2(f_{1}-f_{2})P(X,Y,X,\phi Y)-2f_{4}P(X,Y,hX,\phi Y),
\end{align*}
for any vector fields $X,Y$ orthogonal to $\xi $.
\end{prop}

\section[Contact metric generalized (k,m)-space forms]{Contact metric generalized \protect\boldmath$(\protect\kappa ,
\protect\mu)$-space forms\label{sect-contact}}

In this section we will study contact metric generalized $\left(
\kappa ,\mu \right) $-space forms. The first fundamental fact is
that such a manifold is a generalized $\left( \kappa ,\mu \right)
$-space.

\begin{thm}
\label{th-cmgkmsf-1} If $M(f_{1},\ldots ,f_{6})$ is a contact metric
generalized $\left( \kappa ,\mu \right)$-space form, then it is a
generalized $\left( \kappa ,\mu \right)$-space, with $\kappa=f_{1}-f_{3}$ and $\mu =f_{4}-f_{6}$.
\end{thm}
\begin{proof}
Using the definition of the tensors $R_1,R_2$ and $R_3$, we obtain by direct computation that for every $X,Y$ vector fields on $M$:
\begin{align*}
R_1(X,Y)\xi &=-R_3(X,Y)\xi=\eta(Y)X-\eta(X)Y,\\
R_2(X,Y)\xi &=0.
\end{align*}
Moreover, by the properties \eqref{eq-cont-h} of the tensor $h$ and the definition of $R_4$, $R_5$, $R_6$, it also holds that:
\begin{align*}
R_4(X,Y)\xi &=-R_6(X,Y)\xi=\eta(Y)h X-\eta(X)h Y,\\
R_5(X,Y)\xi &=0.
\end{align*}
Therefore, it would be enough to use the formula \eqref{eq-R-gkmsf} to obtain that the curvature tensor of a generalized $(\kappa,\mu)$-space form satisfies
$$R(X,Y)\xi=(f_1-f_3) \{\eta(Y)X-\eta(X)Y \}+(f_4-f_6) \{\eta(Y)h X-\eta(X)h Y \},$$
for every $X,Y$.
\end{proof}

We know from Theorem~\ref{th-Sas-3} that $M(f_{1},\ldots ,f_{6})$ is
Sasakian if $\kappa =f_{1}-f_{3}=1$. Under the same hypotheses, we
also know that $f_{2}=f_{3}$ (Theorem~\ref{th-Sas-2}) and $h=0$, so
we may take $f_4=f_5=f_6=0$.
Therefore, in the remainder of this section we will study non-Sasakian generalized $(\kappa,\mu)$-space forms $M(f_{1},\ldots ,f_{6})$, that is, those with $\kappa =f_{1}-f_{3}\neq 1$.

We will use the following result from \cite{BKP-95}:

\begin{thm}
\label{th-km-sp-1} If $M$ is a $\left( \kappa ,\mu \right) $-space, then $\kappa \leq 1$. If $\kappa =1$, then $h=0$ and $M$ is a Sasakian
manifold. If $\kappa<1$, then $M$ admits three mutually orthogonal
and integrable distributions ${\mathcal {D}}(0),{\mathcal {D}}(\lambda )$
and ${\mathcal {D}}(-\lambda ) $ determined by the eigenspaces of $h$,
where $\lambda =\sqrt{1-\kappa }$. Moreover,
\begin{align}
R(X_{\lambda },Y_{\lambda })Z_{-\lambda } =&(\kappa -\mu )\{g(\phi
Y_{\lambda },Z_{-\lambda })\phi X_{\lambda }-g(\phi X_{\lambda
},Z_{-\lambda })\phi Y_{\lambda }\},   \label{eq1} \\
R(X_{-\lambda },Y_{-\lambda })Z_{\lambda } =&(\kappa -\mu )\{g(\phi
Y_{-\lambda },Z_{\lambda })\phi X_{-\lambda }-g(\phi X_{-\lambda
},Z_{\lambda })\phi Y_{-\lambda }\},  \\
R(X_{\lambda },Y_{-\lambda })Z_{-\lambda } =&\kappa g(\phi
X_{\lambda },Z_{-\lambda })\phi Y_{-\lambda }+\mu g(\phi X_{\lambda
},Y_{-\lambda
})\phi Z_{-\lambda },   \label{eq3teodim5} \\
R(X_{\lambda },Y_{-\lambda })Z_{\lambda } =&- \kappa g(\phi
Y_{-\lambda },Z_{\lambda })\phi X_{\lambda }-\mu g(\phi Y_{-\lambda
},X_{\lambda
})\phi Z_{\lambda },  \\
R(X_{\lambda },Y_{\lambda })Z_{\lambda } =&(2(1+\lambda )-\mu
)\{g(Y_{\lambda },Z_{\lambda })X_{\lambda }-g(X_{\lambda
},Z_{\lambda
})Y_{\lambda }\},  \\
R(X_{-\lambda },Y_{-\lambda })Z_{-\lambda } =&(2(1-\lambda )-\mu
)\{g(Y_{-\lambda },Z_{-\lambda })X_{-\lambda }-g(X_{-\lambda
},Z_{-\lambda })Y_{-\lambda }\}.  \label{eq2}
\end{align}
\end{thm}

The following results appear in \cite{Kouf-97}:

\begin{thm}
\label{th-km-sp-2} Let $M$ be a $\left( \kappa ,\mu \right) $-space
of dimension greater than or equal to $5$. If the $\phi $-sectional
curvature at any point of $M$ is independent of the choice of the $\phi$-section at that point, then it is constant on $M$ and the curvature
tensor is given by
\begin{equation}
R=\frac{c+3}{4} \, R_{1}+\frac{c-1}{4} \, R_{2}+\left(
\frac{c+3}{4}-\kappa \right) R_{3}+R_{4}+\frac{1}{2} \, R_{5}+\left(
1-\mu \right) R_{6},
\end{equation}
where $c$ is the constant $\phi $-sectional curvature. Moreover, if $\kappa \neq 1$, then $\mu =\kappa +1$ and $c=-\ 2\kappa -1$.
\end{thm}

\begin{thm}
\label{th-km-sp-3} Let $M$ be a non-Sasakian $\left( \kappa ,\mu \right)$-space of dimension greater than or equal to $5$. Then $M$ has constant $\phi $-sectional curvature if and only if $\mu =\kappa +1$.
\end{thm}

Furthermore, in \cite{Kouf-00} the following two theorems have been
proved:

\begin{thm}
\label{th-km-sp-4} Let $M$ be a non-Sasakian generalized $\left(
\kappa ,\mu \right) $-space of dimension greater than or equal to
$5$. Then, the functions $\kappa ,\mu $ are constant, that is, $M$
is a $\left( \kappa ,\mu \right) $-space.
\end{thm}

\begin{thm}\label{th-km-sp-5} Let $M$ be a non-Sasakian generalized $\left(
\kappa ,\mu \right) $-space. If $\kappa ,\mu $ satisfy the condition
$a\kappa +b\mu =c$ (where $a,b,c$ are constant), then $\kappa ,\mu $ are
constant.
\end{thm}

Applying the previous theorems to contact metric generalized $\left(
\kappa ,\mu \right) $-space forms $M(f_{1},\ldots ,f_{6})$, we
deduce:

\begin{thm}\label{th-constant}
If $M(f_{1},\ldots ,f_{6})$ is a non-Sasakian, contact metric generalized $\left( \kappa ,\mu \right) $-space form and $a(f_{1}-f_{3})+b(f_{4}-f_{6})=c$, ($a,b,c$ constant), then $f_{1}-f_{3}$ and $f_{4}-f_{6}$ are
constant.
\end{thm}

We can also obtain the following theorem, in which we prove that the
functions of a contact metric generalized $\left( \kappa ,\mu \right)$-space form $M(f_{1},\ldots ,f_{6})$ of dimension greater than or
equal to 5 are constant and are related. We also obtain a kind of
converse result.

\begin{thm}
\label{th-km-dim5} If $M(f_{1},\ldots ,f_{6})$ is a non-Sasakian,
contact metric generalized $\left( \kappa ,\mu \right) $-space form
of dimension greater than or equal to $5$, then $M$ has constant
$\phi $-sectional curvature $c=2 f_6-1>-3$ and
\begin{align}
f_{1} =&\frac{f_{6}+1}{2}, & f_{2}=&\frac{f_{6}-1}{2}, & f_{3}=&\frac{3f_{6}+1}{2},   \nonumber \\
f_{4} =&1,                 & f_{5}=&\frac{1}{2},       & f_{6}=& \;{\rm constant}>-1, \nonumber \\
\kappa =&-f_{6}<1,          & \mu =&1-f_{6}<2, \label{eq-km-dim5}
\end{align}
Hence $M$ is a $(-f_{6},1-f_{6})$-space with constant $\phi
$-sectional curvature $c=2f_{6}-1>-3$.

Conversely, let $M$ be a $(-f_{6},1-f_{6})$-space of dimension
greater than or equal to $5$ and constant $\phi $-sectional
curvature $c=2f_{6}-1>-3$. Then $M$ is a non-Sasakian, contact
metric generalized $\left( \kappa ,\mu
\right) $-space form $M(f_{1},\ldots ,f_{6})$ with constant functions $f_{1},\ldots ,f_{6}$ satisfying $\left( \ref{eq-km-dim5}\right) $.
\end{thm}

\begin{proof} In view of Theorems~\ref{th-cmgkmsf-1} and \ref{th-km-sp-4}, it follows that $M$ is a generalized $(f_{1}-f_{3},f_{4}-f_{6}) $-space and that $f_{1}-f_{3},f_{4}-f_{6}$
are constant. Therefore, $M$ is a $(f_{1}-f_{3},f_{4}-f_{6})$-space
and applying Theorem~\ref{th-km-sp-1} we obtain that the curvature
tensor satisfies equations (\ref{eq1})--(\ref{eq2}) for $\kappa
=f_{1}-f_{3}$, $\mu =f_{4}-f_{6}$, $X_{\lambda },Y_{\lambda
},Z_{\lambda }\in {\mathcal {D}}(\lambda )$ and $X_{-\lambda
},Y_{-\lambda },Z_{-\lambda }\in {\mathcal {D}}(-\lambda )$.
In particular
\begin{align*}
R(X_{\lambda },Y_{\lambda })Z_{-\lambda }=&(f_{1}-f_{3}-f_{4}+f_{6}) \times\\
&\times \{g(\phi Y_{\lambda },Z_{-\lambda })\phi X_{\lambda }-g(\phi X_{\lambda },Z_{-\lambda })\phi Y_{\lambda}\}, \\
R(X_{\lambda },Y_{-\lambda })Z_{-\lambda } =&(f_{1}-f_{3})g(\phi X_{\lambda },Z_{-\lambda })\phi Y_{-\lambda}\\
& +(f_{4}-f_{6})g(\phi X_{\lambda },Y_{-\lambda })\phi Z_{-\lambda }, \\
R(X_{\lambda },Y_{\lambda })Z_{\lambda }=&(2(1+\sqrt{1-f_{1}+f_{3}})-f_{4} +f_{6}) \times \\
 & \times \{g(Y_{\lambda },Z_{\lambda })X_{\lambda}-g(X_{\lambda},Z_{\lambda })Y_{\lambda }\}, \\
R(X_{-\lambda },Y_{-\lambda })Z_{-\lambda }=&(2(1-\sqrt{1-f_{1}+f_{3}})-f_{4} +f_{6}) \times \\
  &\times \{g(Y_{-\lambda},Z_{-\lambda})X_{-\lambda}-g(X_{-\lambda },Z_{-\lambda })Y_{-\lambda }\}.
\end{align*}

On the other hand, using the definition of generalized $\left(
\kappa ,\mu \right) $-space form and the properties of contact
metric manifolds we get:
\begin{align*}
R(X_{\lambda },Y_{\lambda })Z_{-\lambda } =&(f_{2}-f_{5}(1-f_{1}+f_{3})) \times\\
& \times \{g(\phi Y_{\lambda },Z_{-\lambda})\phi X_{\lambda }
-g(\phi X_{\lambda },Z_{-\lambda })\phi Y_{\lambda }\}, \\
R(X_{\lambda },Y_{-\lambda })Z_{-\lambda }
=&-(f_{2}+f_{5}(1-f_{1}+f_{3}))g(\phi X_{\lambda },Z_{-\lambda})\phi Y_{-\lambda } \\
&-2  f_{2}g(\phi X_{\lambda },Y_{-\lambda })\phi Z_{-\lambda}\\
&  +(f_{1}-f_{5}(1-f_{1}+f_{3}))g(Y_{-\lambda},Z_{-\lambda})X_{\lambda}, \\
R(X_{\lambda },Y_{\lambda })Z_{\lambda}
=&\left(f_{1}+f_{5}(1-f_{1}+f_{3})+2f_{4}\sqrt{1-f_{1}+f_{3}} \right) \times\\
&\hfill \times \{g(Y_{\lambda },Z_{\lambda })X_{\lambda }-g(X_{\lambda},Z_{\lambda})Y_{\lambda }\}, \\
R(X_{-\lambda },Y_{-\lambda })Z_{-\lambda }
=&\left(f_{1}+f_{5}(1-f_{1}+f_{3})-2f_{4}\sqrt{1-f_{1}+f_{3}} \right) \times \\
& \hfill \times  \{g(Y_{-\lambda },Z_{-\lambda })X_{-\lambda }-g(X_{-\lambda},Z_{-\lambda })Y_{-\lambda }\}.
\end{align*}

Combining both sets of equations we can write:
\begin{align}
&(f_{1}-f_{3}-f_{4}+f_{6})\{g(\phi Y_{\lambda },Z_{-\lambda })\phi
X_{\lambda }-g(\phi X_{\lambda },Z_{-\lambda })\phi Y_{\lambda }\}\nonumber \\
&\quad =(f_{2}-f_{5}(1-f_{1}+f_{3}))\{g(\phi Y_{\lambda
},Z_{-\lambda })\phi X_{\lambda }-g(\phi X_{\lambda },Z_{-\lambda })\phi Y_{\lambda }\},  \label{eq-1}
\end{align}
\begin{align}
(&f_{1}-f_{3})g(\phi X_{\lambda },Z_{-\lambda })\phi Y_{-\lambda }+(f_{4}-f_{6})g(\phi X_{\lambda },Y_{-\lambda })\phi Z_{-\lambda } \nonumber \\
& =-\ (f_{2}+f_{5}(1-f_{1}+f_{3}))g(\phi X_{\lambda },Z_{-\lambda })\phi Y_{-\lambda }-2f_{2}g(\phi X_{\lambda },Y_{-\lambda })\phi Z_{-\lambda }   \nonumber \\
&\quad +\ (f_{1}-f_{5}(1-f_{1}+f_{3}))g(Y_{-\lambda },Z_{-\lambda })X_{\lambda }, \label{eq-3}
\end{align}
\begin{align}
(2(1+&\sqrt{1-f_{1}+f_{3}})-f_{4}+f_{6})\{g(Y_{\lambda },Z_{\lambda
})X_{\lambda }-g(X_{\lambda },Z_{\lambda })Y_{\lambda }\} \nonumber\\
& =\left( f_{1}+f_{5}(1-f_{1}+f_{3})+2f_{4}\sqrt{1-f_{1}+f_{3}} \right) \times \nonumber\\
&\times  \{g(Y_{\lambda},Z_{\lambda})X_{\lambda} -g(X_{\lambda},Z_{\lambda })Y_{\lambda }\},\label{eq-5}
\end{align}
\begin{align}
(2(1- & \sqrt{1-f_{1}+f_{3}})-f_{4}+f_{6})\{g(Y_{-\lambda },Z_{-\lambda
})X_{-\lambda }-g(X_{-\lambda },Z_{-\lambda })Y_{-\lambda }\} \nonumber\\
&=\left( f_{1}+f_{5}\left( 1-f_{1}+f_{3}\right) -2f_{4}\sqrt{1-f_{1}+f_{3}}\right) \times \nonumber\\
&  \times \hfill \left\{ g\left( Y_{-\lambda },Z_{-\lambda }\right)X_{-\lambda } -g\left( X_{-\lambda},Z_{-\lambda}\right) Y_{-\lambda} \right\}. \label{eq-6}
\end{align}

The dimension is greater than or equal to $5$, so we can take two
mutually
orthogonal unit vector fields $X_{\lambda },Y_{\lambda }\in {\mathcal {D}} (\lambda )$ in the equation (\ref{eq-1}) and choosing  $Z_{-\lambda }=\phi X_{\lambda }$ we obtain
\[
(f_{1}-f_{3}-f_{4}+f_{6})(-\phi Y_{\lambda
})=(f_{2}-f_{5}(1-f_{1}+f_{3}))(-\phi Y_{\lambda }),
\]
and we deduce that
\begin{equation}
f_{1}-f_{3}-f_{4}+f_{6}=f_{2}-f_{5}(1-f_{1}+f_{3}).
\label{{eq-1.1}}
\end{equation}

From (\ref{eq-3}) we obtain three equations, depending on the choice
we make. If we take two mutually orthogonal unit vector fields
$Y_{-\lambda },Z_{-\lambda }\in {\mathcal {D}}(-\lambda )$ and pick
$X_{\lambda }=\phi Z_{-\lambda }$, we obtain:
\begin{equation}
-(f_{1}-f_{3})=f_{2}+f_{5}(1-f_{1}+f_{3}).  \label{eq-3.1}
\end{equation}
For $Y_{-\lambda },Z_{-\lambda }\in {\mathcal {D}}(-\lambda )$
orthogonal and unit and $X_{\lambda }=\phi Y_{-\lambda }$, it
follows that
\begin{equation}
-(f_{4}-f_{6})=2f_{2}.  \label{eq-3.2}
\end{equation}
If we now take $Y_{-\lambda }=Z_{-\lambda }$ unit and $X_{\lambda
}=\phi Y_{-\lambda }$, we have the following equation:
\begin{equation}
-(f_{1}-f_{3})-(f_{4}-f_{6})=f_{1}+3f_{2}.  \label{eq-3.3}
\end{equation}

In (\ref{eq-5}) we choose two unit vector fields $Y_{\lambda
} \bot Z_{\lambda }\in {\mathcal {D}}(\lambda )$ and $X_{\lambda }=Z_{\lambda
}$ and we obtain
\begin{equation}
2(1+\sqrt{1-f_{1}+f_{3}})-f_{4}+f_{6}=f_{1}+f_{5}(1-f_{1}+f_{3})+2f_{4}\sqrt{1-f_{1}+f_{3}}.  \label{eq-5.1}
\end{equation}

Last, if we take two unit vector fields
$X_{-\lambda} \bot Z_{-\lambda}\in {\mathcal {D}}(-\lambda)$ and
$Y_{-\lambda}=Z_{-\lambda}$ in (\ref{eq-6}) then
\begin{equation}
2 ( 1 - \sqrt{1-f_1+f_3})-f_4+ f_6 =
f_1+f_5(1-f_1+f_3)-2f_4\sqrt{1-f_1+f_3}. \label{eq-6.1}
\end{equation}

Combining and reordering equations (\ref{{eq-1.1}})-(\ref{eq-6.1}),
we get the following compatible system
\[
\left.
\begin{array}{rcc}
f_{1}-f_{2}-f_{3}-f_{4}+f_{6}+f_{5}(1-f_{1}+f_{3}) & = & 0 \\
f_{1}+f_{2}-f_{3}+f_{5}(1-f_{1}+f_{3}) & = & 0 \\
2f_{2}+f_{4}-f_{6} & = & 0 \\
2f_{1}+3f_{2}-f_{3}+f_{4}-f_{6} & = & 0 \\
f_{5}(1-f_{1}+f_{3})+2(f_{4}-1)\sqrt{1-f_{1}+f_{3}}+f_{1}+f_{4}-f_{6}-2
& = & 0 \\
f_{5}(1-f_{1}+f_{3})+2(1-f_{4})\sqrt{1-f_{1}+f_{3}}+f_{1}+f_{4}-f_{6}-2
& = & 0
\end{array}
\right\}
\]
whose solution is:
\begin{align*}
f_{1}&=\displaystyle\frac{f_{6}+1}{2}, & f_{2}&=\displaystyle\frac{f_{6}-1}{2}, & f_{3}&=\displaystyle\frac{3f_{6}+1}{2}, \\
f_{4}&=1, &    f_{5}&=\displaystyle\frac{1}{2},&  f_{6} &{\rm  \quad arbitrary}.
\end{align*}

Therefore, $\kappa =f_{1}-f_{3}=-f_{6}$, $\mu =f_{4}-f_{6}=1-f_{6}$ and $c=f_{1}+3f_{2}=2f_{6}-1$, by virtue of
Proposition~\ref{prop-gkmsf-K(X,phiX)} and
Theorem~\ref{th-cmgkmsf-1}.

Now, since $\kappa $ is constant and less than $1$, we have that
$f_{6}$ is constant and greater than $-1$ and we achieve the result.

Conversely, given a $(-f_{6},1-f_{6})$-space of dimension greater than or equal to 5 with $\phi $-sectional curvature $c=2f_{6}-1>-3$, it follows that $f_{6}$
must
be a constant function greater than $-1$. It is sufficient to apply Theorem~\ref{th-km-sp-2} with $c=2f_{6}-1$ to get that the manifold has
curvature tensor
\begin{equation}
R=\frac{f_{6}+1}{2}R_{1}+\frac{f_{6}-1}{2}R_{2}+\frac{3f_{6}+1}{2} R_{3}+R_{4}+\frac{1}{2}R_{5}+f_{6}R_{6},
\end{equation}
so that it is a contact metric generalized $\left( \kappa ,\mu \right)$-space form $M(f_{1},\ldots ,f_{6})$ with functions $f_{1},\ldots,f_{6}$ satisfying equations (\ref{eq-km-dim5}). It is obvious that
the manifold is non-Sasakian because $\kappa =f_{1}-f_{3}=-f_{6}<1$.
\end{proof}

\begin{rem-new}
We observe that $f_{4},f_{5}\neq 0$ in (\ref{eq-km-dim5}), so there
are no examples of non-Sasakian, contact metric generalized Sasakian
space forms of dimension greater than or equal to $5$ (already seen
in \cite{AC-DGA}).\
\end{rem-new}

We will now give a method to construct $(-f_{6},1-f_{6})$-spaces
with
constant $\phi $-sectional curvature $c=2f_{6}-1$ for every constant $f_{6}>-1$. Due to the previous theorem, they will be examples of
contact metric generalized $\left( \kappa ,\mu \right) $-space
forms.

Let $M$ be a manifold of dimension greater than or equal to 5 and
constant sectional curvature $c_{s}>-1(c_{s}\neq 1)$. Then its
tangent sphere bundle with the usual contact metric, $(T_{1}M,\xi_1, \eta_1 ,\phi_1 ,g_1)$, is a $\left( \kappa ,\mu
\right) $-space with $\kappa =c_{s}(2-c_{s})\neq 1$ and $\mu =-2c_{s}<2$ (\cite[Theorem~4]{BKP-95}).

By applying a $D_{a}$-homothetic deformation we obtain $(T_{1}M,\overline{\xi},\overline{\eta},\overline{\phi},\overline{g})$, with
\[
\overline{\xi }=\frac{1}{a}\,\xi ,\quad \overline{\eta }=a\eta
,\quad \overline{\phi }=\phi \quad {\rm and}\quad
\overline{g}=ag+a(a-1)\eta \otimes \eta,
\]
where $a>0$ is a real number.

It is known from \cite{BKP-95} that this is a $(\overline{\kappa
},\overline{\mu })$-space with
\[
\overline{\kappa }=\frac{\kappa +a^{2}-1}{a^{2}}\neq 1\quad {\rm
and}\quad \overline{\mu} = \frac{\mu +2a-2}{a}.
\]

If we choose $a= (\kappa -1)/(\mu -2)>0$, then $\overline{\mu }=\overline{\kappa }+1$ and $(T_{1}M,\overline{\xi },\overline{\eta },\overline{\phi },\overline{g})$ has constant $\phi $-sectional curvature $\overline{c}=-(\overline{\kappa }+\overline{\mu })$ because of
Theorems~\ref{th-km-sp-2} and \ref{th-km-sp-3}.

Therefore, $\overline{\kappa }=-f_{6}$, $\overline{\mu }=1-f_{6}$ and $\overline{c}=2f_{6}-1$ if and only if
\begin{equation}
\left( 3-f_{6}\right) c_{s}^{2}+\left( 10+2f_{6}\right) c_{s}+\left(
3-f_{6}\right) =0.  \label{eq-sect-curv}
\end{equation}
If $f_{6}=3$, equation (\ref{eq-sect-curv}) has solution $c_{s}=0$,
which is
in particular greater than $-1$ and not equal to 1. If $f_{6}\neq 3$, then (\ref{eq-sect-curv}) has the real solutions
\begin{equation}
c_{s}=\frac{-5-f_{6}\pm 4\sqrt{f_{6}+1}}{3-f_{6}},  \label{eq-cs}
\end{equation}
which are not equal to 1 because $f_{6}$ is greater than $-1$.
Furthermore, if we consider the positive sign on (\ref{eq-cs}), it
can be proved that $c_{s}$ is also greater than $-1$. Therefore, we
have managed to obtain examples of $(-f_{6},1-f_{6})$-spaces with
the required conditions for every constant function $f_{6}>-1$.

\begin{rem-new}
\label{rem-alt-proof} An alternative proof of
Theorem~\ref{th-km-dim5} follows from \cite[Theorem~5]{Boeckx-99},
which states that, if $M$ is a
non-Sasakian $\left( \kappa ,\mu \right) $-space, then its curvature tensor $R$ is given by
\begin{equation}\label{eq-R-non-Sas-km}
\begin{aligned}
R( X,Y, Z, W) &=\left( 1-\frac{\mu }{2}\right) R_{1}(X,Y,Z,W) -\frac{\mu }{2}\,R_{2}( X,Y,Z,W)  \\
&+\left( 1-\frac{\mu }{2}-\kappa \right) R_{3}( X,Y,Z,W)+R_{4}( X,Y,Z,W)  \\
&+\left( \frac{1-\frac{\mu }{2}}{1-\kappa }\right)
 \left\{ g(hY,Z) g( h X,W) -g( hX,Z) g(hY,W) \right\}  \\
&+\left( \frac{\kappa-\frac{\mu }{2} }{1-\kappa }\right)
\left\{g( \phi hY,Z)g( \phi hX,W) -g( \phi hX,Z) g( \phi hY,W) \right\}   \\
&+\left( 1-\mu \right) R_{6}( X,Y,Z,W).
\end{aligned}
\end{equation}
Therefore, if $M(f_{1},\ldots ,f_{6})$ is a contact metric generalized $\left( \kappa ,\mu \right) $-space form of dimension greater than or
equal to $5$ satisfying $f_{1}-f_{3}<1$, then $M$ is a non-Sasakian
$\left( \kappa ,\mu \right) $-space with $\kappa =f_{1}-f_{3}<1$ and
$\mu =f_{4}-f_{6}$
thanks to Theorems~\ref{th-cmgkmsf-1} and \ref{th-km-sp-4}. Comparing (\ref{eq-R-gkmsf}) with (\ref{eq-R-non-Sas-km}), we obtain a system whose
solution is given by (\ref{eq-km-dim5}).
\end{rem-new}

\

What can we say now for $3$-dimensional generalized $\left( \kappa,\mu \right) $-space form $M^{3}(f_{1},\ldots ,f_{6})$? First, let us mention that the writing of its curvature tensor is
not unique:

\begin{thm}
\label{th-R-non-unique-3d} Let $M^{3}$ be a contact metric manifold
such that its curvature tensor can be simultaneously written as
\begin{equation}
R=f_{1}R_{1}+f_{2}R_{2}+f_{3}R_{3}+f_{4}R_{4}+f_{5}R_{5}+f_{6}R_{6}
\label{eq-R-non-unique-3d-1}
\end{equation}
and
\begin{equation}
R=f_{1}^{\ast }R_{1}+f_{2}^{\ast }R_{2}+f_{3}^{\ast
}R_{3}+f_{4}^{\ast }R_{4}+f_{5}^{\ast }R_{5}+f_{6}^{\ast }R_{6},
\label{eq-R-non-unique-3d-2}
\end{equation}
where $f_{1}-f_{3}<1$. Then the functions $f_{i}$ and $f_{i}^{\ast }$, $i=1,\ldots ,6$, are related as follows,
\begin{equation}
\left.
\begin{array}{c}
f_{1}^{\ast }=f_{1}+f,\quad f_{2}^{\ast }=f_{2}-f/3,\quad
f_{3}^{\ast
}=f_{3}+f,  \\
f_{4}^{\ast }=f_{4}+\overline{f},\quad f_{6}^{\ast }=f_{6}+\overline{f},
\end{array}
\right.  \label{eq-R-non-unique-3d-3}
\end{equation}
where $f$ and $\overline{f}$ are arbitrary functions on $M$.
\end{thm}

\begin{proof} We know that the manifold is in particular a generalized $(\kappa,\mu)$-space with $\kappa=f_1-f_3<1$ and $\mu=f_4-f_6$ (Theorem \ref{th-cmgkmsf-1}). Therefore, we can consider a $\phi $-basis $\{X,\phi X,\xi \}$ with $X \in D(\lambda)$, where $\lambda=\sqrt{1-\kappa}=\sqrt{1-(f_1-f_3)}>0$ (Theorem \ref{th-km-sp-1}). By virtue of Propositions \ref{prop-gkmsf-K(X,phiX)}, \ref{prop-gkmsf-K(X,xi)} and Corollary \ref{cor-gkmsf-K(phiX,xi)}, if we calculate $K(X,\xi)$, $K(\phi X,\xi)$ and $K(X,\phi X)$ by using both \eqref{eq-R-non-unique-3d-1} and \eqref{eq-R-non-unique-3d-2} we obtain the system
\[
\left\{
\begin{array}{rcl}
(f_{1}^{\ast }-f_{1})-(f_{3}^{\ast }-f_{3}) & = & 0 \\
(f_{4}^{\ast }-f_{4})-(f_{6}^{\ast }-f_{6}) & = & 0 \\
(f_{1}^{\ast }-f_{1})+3(f_{2}^{\ast }-f_{2}) & = & 0
\end{array}
\right.
\]
whose general solution is given by (\ref{eq-R-non-unique-3d-3}).
\end{proof}

\begin{rem-new}
In the conditions of the previous theorem, if $\kappa =f_{1}-f_{3}=1$, then $M^{3}$ is a Sasakian manifold and therefore it is a generalized
Sasakian space form $M(f_{1},f_{2},f_{3})$. In \cite{AC-DGA}, P.
Alegre and A.
Carriazo proved that in such a case the functions $f_{i}$ and $f_{i}^{\ast }$, $i=1,2,3$, are related as in (\ref{eq-R-non-unique-3d-3}).
\end{rem-new}

The converse of Theorem~\ref{th-R-non-unique-3d} is also true:

\begin{thm}
Let $M^{3}(f_{1},\ldots ,f_{6})$ be a contact metric generalized
$\left( \kappa ,\mu \right) $-space form. If we define the functions
$f_{1}^{\ast
},f_{2}^{\ast },f_{3}^{\ast },f_{4}^{\ast },f_{6}^{\ast }$ as in (\ref{eq-R-non-unique-3d-3}), for certain functions $f,\overline{f}$ on
$M$, and
we take an arbitrary function $f_{5}^{\ast }$, then $M^{3}$ is also a generalized $\left( \kappa ,\mu \right) $-space form $M^{3}(f_{1}^{\ast },\ldots
,f_{6}^{\ast })$.
\end{thm}

\begin{proof} It follows from (\ref{eq-R-gkmsf}) and (\ref{eq-R-non-unique-3d-3}) that the curvature tensor satisfies
\[
R=\sum_{i=1}^{6}f_{i}R_{i}=\sum_{i=1}^{6}f_{i}^{\ast }R_{i}+
f\left( -R_{1}+\frac{1}{3}R_{2}-R_{3}\right) -\overline{f}(
R_{4}+R_{6}) +( f_{5}-f_{5}^{\ast }) R_{5}.
\]
To obtain (\ref{eq-R-non-unique-3d-2}) it is enough to check that
the last
terms vanish, which is true because $$-R_{1}+\frac{1}{3} R_{2}-R_{3}=R_{4}+R_{6}=R_{5}=0$$ for every $3$-dimensional contact
metric manifold.
\end{proof}

Therefore, if $M^{3}(f_{1},\ldots ,f_{6})$ is a contact metric generalized $\left( \kappa ,\mu \right) $-space form, its curvature tensor can be
written as
\begin{equation}
R=f_{1}^{\ast }R_{1}+f_{3}^{\ast }R_{3}+f_{4}^{\ast }R_{4},
\label{eq-R-cmgkmsf-3d-1}
\end{equation}
so $M^{3}$ is also $M(f_{1}^{\ast },0,f_{3}^{\ast },f_{4}^{\ast },0,0)$ for $f_{1}^{\ast }=f_{1}+3f_{2},f_{3}^{\ast }=f_{3}+3f_{2},f_{4}^{\ast
}=f_{4}-f_{6}.$ In order to consider a unique writing of the
curvature tensor of a 3-dimensional, contact metric generalized
$\left( \kappa ,\mu \right) $-space form, we will choose $R$
satisfying $f_{2}^{\ast }=f_{5}^{\ast }=f_{6}^{\ast }=0$.

D. E. Blair, T. Koufogiorgos and B. J. Papantoniou classified the
3-dimensional $\left( \kappa ,\mu \right) $-spaces in
\cite[Theorem~3]{BKP-95}. Using that result and
Theorems~\ref{th-cmgkmsf-1} and \ref{th-km-sp-1} we get:

\begin{thm}
\label{th-caso3} Let $M^{3}(f_{1},\ldots ,f_{6})$ be a contact
metric generalized $\left( \kappa ,\mu \right) $-space form. If
$f_{1}-f_{3}\neq 1$ is constant, then  $f_{4}-f_{6}$ is also constant, $M$
satisfies
\begin{equation}
2f_{1}+3f_{2}-f_{3}+f_{4}-f_{6}=0,  \label{eq-caso3}
\end{equation}
and it is locally isometric to one of the following Lie groups with
a left invariant metric: $SU(2)(or$ $SO(3))$, $SL(2,{\mathbf R})(or$
$O(1,2))$, $E(2)$ (the group of rigid motions of the Euclidean
$2$-space) or $E(1,1)$ (the group of rigid motions of the Minkowski
$2$-space).

Moreover this structure can occur on:

\begin{itemize}
\item $SU(2)$ or $SO(3)$ if $1-\lambda-\frac{\mu}{2} > 0$ and $1+\lambda-\frac{\mu}{2}> 0$,

\item $SL(2,{\mathbf R})$ or $O(1,2)$ if $\left\{
\begin{array}{llll}
\quad & 1-\lambda -\frac{\mu }{2}<0 & and & 1+\lambda -\frac{\mu }{2}>0 \\
{\rm or} &  &  &  \\
\quad & 1-\lambda -\frac{\mu }{2}<0 & and & 1+\lambda -\frac{\mu }{2}<0,
\end{array}
\right. $

\item $E(2)$ if $1-\lambda-\frac{\mu}{2} =0$ and $\mu < 2$,

\item $E(1,1)$ if $1+\lambda-\frac{\mu}{2} =0$ and $\mu > 2$,
\end{itemize}
where $\kappa=f_1-f_3, \lambda=\sqrt{1-\kappa}$ and
$\mu=-2f_1-3f_2+f_3$.
\end{thm}

\begin{proof} Theorems \ref{th-Sas-2} and \ref{th-cmgkmsf-1}  imply that $M$ is a non-Sasakian generalized $\left( \kappa ,\mu \right) $-space with $\kappa =f_{1}-f_{3}$ and $\mu =f_{4}-f_{6}$. Since $f_1-f_3$ is constant, then we know from Theorem \ref{th-constant} that  $f_{4}-f_{6}$ is also constant and hence $M$ is a $(\kappa,\mu)$-space.  Applying Theorem~\ref{th-km-sp-1} we obtain
that $R$ must satisfy the six equations (\ref{eq1})-(\ref{eq2}).

If we take unit vector fields $Y_{-\lambda }=Z_{-\lambda }$ and
$X_{\lambda }=\phi Y_{-\lambda }$ in (\ref{eq3teodim5}), then
\[
-\ (f_{1}-f_{3})-(f_{4}-f_{6})=f_{1}+3f_{2}.
\]

A short calculation yields $2f_{1}+3f_{2}-f_{3}+f_{4}-f_{6}=0$ and
$\mu =-2f_{1}-3f_{2}+f_{3}$. We only have to apply
\cite[Theorem~3]{BKP-95} to end the proof. \end{proof}

\begin{rem-new}
We notice that the different cases of Theorem~\ref{th-caso3} depend
on the
value of $\kappa $ and $\mu $, which are determined only by the functions $f_{1},f_{2}$ and $f_{3}$ and do not depend explicitly on the functions
$f_{4},f_{5}$ or $f_{6}$.
\end{rem-new}

Let us recall that a contact metric manifold is said to be {\em $\eta$-Einstein} \cite[p.~105]{Blair-02} if it satisfies
\[
Q=aI+b\eta \otimes \xi ,
\]
where $a,b$ are some differentiable functions on $M$.

In \cite{BKS-90}, D. E. Blair, T. Koufogiorgos and R. Sharma studied
the contact metric manifolds satisfying $Q\phi =\phi Q$ and obtained
the following result:

\begin{prop}
If $( M^{3},\phi ,\xi ,\eta ,g) $ is a contact metric
manifold, then the following statements are equivalent:

\begin{itemize}
\item[{\rm (i)}] $M$ is $\eta $-Einstein,

\item[{\rm (ii)}] $Q\phi =\phi Q$, where $Q$ is the Ricci
operator,

\item[{\rm (iii)}] $M$ is a $(\kappa ,0)$-space, with $\kappa $
constant.
\end{itemize}
\end{prop}

In view of Theorem~\ref{th-cmgkmsf-1}, we know that a
generalized $\left( \kappa ,\mu \right)$-space form $M(f_{1},\ldots, f_{6})$ with contact metric structure
satisfies condition (iii) if and only if $f_{1}-f_{3}$ is constant and $f_{4}-f_{6}=0$. We will now study when the conditions (i) and (ii)
hold. First, we calculate the Ricci operator $Q$:

\begin{prop}
\label{prop-Q} If $M^{2n+1}(f_{1},\ldots ,f_{6})$ is a contact
metric generalized $\left( \kappa ,\mu \right) $-space form, then
\begin{equation}
Q=\left( 2nf_{1}+3f_{2}-f_{3}\right) I-\left( 3f_{2}+\left(2n-1\right) f_{3} \right) \eta \otimes \xi
+\left( \left( 2n-1\right) f_{4} -f_{6}\right) h. \label{eq-Q-1}
\end{equation}
Moreover, if we also suppose that $\kappa =f_{1}-f_{3}\neq 1$ is constant, then $M$ is a non-Sasakian $\left(
\kappa ,\mu \right) $-space and
\begin{equation}
Q=\left( 2\left( n-1\right) -n\mu \right) I+\left( 2\left(
1-n\right) +n\left( 2\kappa +\mu \right) \right) \eta \otimes \xi
+\left( 2\left( n-1\right) +\mu \right) h.  \label{eq-Q-2}
\end{equation}
\end{prop}

\begin{proof} A straightforward computation with respect to a $\phi $-basis gives (\ref{eq-Q-1}). On the other hand, if $\kappa =f_{1}-f_{3}\neq 1$ is constant, then  $\mu =f_{4}-f_{6}$ is also constant by Theorem \ref{th-constant}. Therefore, the equations (\ref{eq-km-dim5}) hold if $M$ has dimension greater than or equal to 5 and (\ref{eq-caso3}) is true if $M$ is a $3$-dimensional manifold.
Applying them to (\ref{eq-Q-1}) yields (\ref{eq-Q-2}). \end{proof}

\begin{rem-new}
Let us notice that the previous proposition means that a contact
metric generalized $\left( \kappa ,\mu \right) $-space form
$M^{2n+1}(f_{1},\ldots ,f_{6})$ is $\eta $-Einstein if and only if
$f_{4}(2n-1)-f_{6}=0$. In particular, $M^{3}(f_1,\ldots,f_6)$ satisfies condition
(i) if and only if $f_{4}-f_{6}=0$.
\end{rem-new}

We can also prove the following:

\begin{prop}
\label{prop-Qphi-phiQ} If $M^{2n+1}(f_{1},\ldots ,f_{6})$ is a
contact metric generalized $\left( \kappa ,\mu \right) $-space form,
then
\begin{equation}
Q\phi -\phi Q=2((2n-1)f_{4}-f_{6})h\phi , \label{eq-Qphi-phiQ}
\end{equation}
where $Q$ denotes the Ricci operator on $M$.
\end{prop}

\begin{proof} Using $\eta \circ \phi =0$, from
(\ref{eq-Q-1}) we obtain
\begin{equation}
Q\phi =\left( 2nf_{1}+3f_{2}-f_{3}\right) \phi +\left( \left(
2n-1\right) f_{4}-f_{6}\right) h\phi .  \label{eq-Qphi}
\end{equation}
Applying $\phi $ to (\ref{eq-Q-1}) and using (\ref{eq-cont-h}) we
get
\begin{equation}
\phi Q=\left( 2nf_{1}+3f_{2}-f_{3}\right) \phi -\left( \left(
2n-1\right)f_{4} -f_{6}\right) h\phi .  \label{eq-phiQ}
\end{equation}
Therefore, (\ref{eq-Qphi}) and (\ref{eq-phiQ}) imply
(\ref{eq-Qphi-phiQ}). \end{proof}
\begin{rem-new}
We deduce from Proposition~\ref{prop-Qphi-phiQ} that if $M^{2n+1}(f_{1},\ldots ,f_{6})$ is a contact metric generalized $\left( \kappa ,\mu \right)$-space form, then $Q\phi =\phi Q$ is true if and only if $(2n-1)f_{4}-f_{6}=0$. In particular, $M^{3}$ satisfies condition
$Q\phi =\phi Q$  if and only if $f_{4}-f_{6}=0$.
\end{rem-new}

We may resume the previous results in this proposition:

\begin{prop}
\label{prop-clasification-dim3} If $M^{3}(f_{1},\ldots ,f_{6})$ is a
contact metric generalized $\left( \kappa ,\mu \right) $-space form,
then the following conditions are equivalent:

\begin{enumerate}
\item[{\rm (i)}] $M^{3}$ is $\eta $-Einstein,

\item[{\rm (ii)}] $Q\phi =\phi Q$, where $Q$ denotes the Ricci
operator,

\item[{\rm (iii)}] $M^{3}$ is a $\left( f_{1}-f_{3},0\right) $-space,

\item[{\rm (iv)}] $f_{4}-f_{6}=0$.
\end{enumerate}
\end{prop}

\begin{rem-new}
If a contact metric generalized $\left( \kappa ,\mu \right) $-space $M^{3}(f_{1},\ldots ,f_{6})$ satisfies $f_{4}-f_{6}=0$, then
$f_{1}-f_{3}$ is
constant. This would be a particular case of {\rm \cite[Theorem~10]{Sharma-95}}.
\end{rem-new}

D. E. Blair and H. Chen proved in \cite{BC-92} the following
theorem, which improves the classification of contact metric
manifolds satisfying $Q \phi =\phi Q$ given in \cite{BKS-90}:

\begin{thm}
\label{th-clasification-dim3} Let $M^3$ be a contact metric manifold
satisfying $Q\phi=\phi Q$. Then $M^3$ is either Sasakian, flat or
locally isometric to a left invariant metric on the Lie group $SU(2)$ or $SL(2,{\mathbf R})$. In the latter case $M^3$ has constant
$\xi$-sectional curvature $\kappa<1$ and constant $\phi$-sectional
curvature $-\kappa$
(this structure can occur on $SU(2)$ if $\kappa>0$ and on $SL(2,{\mathbf R})$ if $\kappa <0)$.
\end{thm}

We deduce from the previous theorem the following proposition:

\begin{prop}
Let $M^{3}(f_{1},\ldots ,f_{6})$ be a contact metric generalized
$\left(
\kappa ,\mu \right) $-space form. If $f_{1}-f_{3}\neq 1$ and $f_{4}-f_{6}=0$, then $f_{1}-f_{3}$ and $f_{1}+3f_{2}$ are constant and $2f_{1}+3f_{2}-f_{3}=0$ holds.
\end{prop}

\begin{proof}
We know that $M^3$ is a contact metric
manifold
satisfying $Q\phi=\phi Q$ due to Proposition~\ref{prop-clasification-dim3}. Applying
Theorem~\ref{th-clasification-dim3} and
using the $\phi$-sectional and $\xi$-sectional formulas from Propositions~\ref{prop-gkmsf-K(X,phiX)} and \ref{prop-gkmsf-K(X,xi)} we get the
wanted result.
\end{proof}

We will now study the scalar curvature:

\begin{prop}
\label{prop-tau} If $M^{3}(f_{1},\ldots ,f_{6})$ is a contact metric
generalized $\left( \kappa ,\mu \right) $-space form, then the
scalar curvature $\tau $ is given by
\begin{equation}
\tau =2(3f_{1}+3f_{2}-2f_{3}).  \label{eq-tau}
\end{equation}
Moreover, if $\kappa =f_{1}-f_{3}\neq 1$ is constant, then
\[
\tau =2(\kappa -\mu).
\]
\end{prop}

\begin{proof} A straightforward computation with respect to a $\phi $-basis yields
\[
\tau =tr (Q)=2 \left( K\left( X,\phi X\right) +K\left( X,\xi \right) +K\left( \phi X,\xi \right) \right) ,
\]
and using Propositions~\ref{prop-gkmsf-K(X,phiX)},
\ref{prop-gkmsf-K(X,xi)} and Corollary~\ref{cor-gkmsf-K(phiX,xi)} we
get (\ref{eq-tau}).

Furthermore, if $\kappa =f_{1}-f_{3}\neq 1$ is constant, then   Theorem \ref{th-constant} implies that $\mu =f_{4}-f_{6}$ is also constant  and we only need to apply (\ref{eq-caso3}) to obtain $3f_{1}+3f_{2}-2f_{3}=(f_{1}-f_{3})-(f_{4}-f_{6})=\kappa -\mu .$
 \end{proof}

\begin{cor}
Let $M^{3}(f_{1},\ldots ,f_{6})$ be a contact metric generalized
$\left( \kappa ,\mu \right) $-space form with $\kappa
=f_{1}-f_{3}\neq 1$ constant. Then
\[
Q=\left( \frac{\tau}{2} -\kappa \right) I+\left( 3 \kappa -\frac{\tau}{2} \right) \eta
\otimes \xi +\left( \kappa -\frac{\tau}{2} \right) h.
\]
\end{cor}

We will prove the next theorem using the formulas we have obtained
for $Q$ and $\tau $ and the expression that connects both of them to
$R$ on a 3-dimensional Riemannian manifold:
\begin{align}
R\left( X,Y\right) Z =&g\left( Y,Z\right) QX-g\left( X,Z\right)
QY+g(QY,Z)X-g(QX,Z)Y  \nonumber \\
&-\frac{\tau}{2} \{g\left( Y,Z\right) X-g\left( X,Z\right) Y\}
\label{eq-R-3d}
\end{align}
for any vector fields $X,Y,Z$ on $M$.

\begin{thm}
If $M^{3}(f_{1},\ldots ,f_{6})$ is a contact metric generalized
$\left( \kappa ,\mu \right) $-space form, then its curvature tensor
can be written as:
\begin{equation}
R=(f_{1}+3f_{2})R_{1}+(3f_{2}+f_{3})R_{3}+(f_{4}-f_{6})R_{4}.  \label{solounavez}
\end{equation}
If we also suppose that $\kappa =f_{1}-f_{3}\neq 1$ is constant, then
\[
R=-(\kappa +\mu )R_{1}-(2\kappa +\mu )R_{3}+\mu R_{4}.
\]
\end{thm}

\begin{proof} Equation (\ref{solounavez}) is obtained by
substituting (\ref{eq-Q-1}) and (\ref{eq-tau}) in (\ref{eq-R-3d}).
Moreover,
if $\kappa =f_{1}-f_{3}\neq 1$ is constant, then we obtain from Thereom \ref{th-constant} that $\mu =f_{4}-f_{6}$ is also constant. Therefore, (\ref{eq-caso3}) holds, which yields $f_{1}+3f_{2}=-(\kappa +\mu )$ and $3f_{2}+f_{3}=-(2\kappa -\mu )$. \end{proof}

We will now check that the first example of a generalized $\left(
\kappa ,\mu \right) $-space of dimension 3 given by T. Koufogiorgos
and C. Tsichlias in
\cite{Kouf-00} is a contact metric generalized $\left( \kappa ,\mu \right)$-space form $M^{3}(f_{1}^{\ast },0,f_{3}^{\ast },f_{4}^{\ast },0,0)$ with $f_{1}^{\ast },f_{3}^{\ast },f_{4}^{\ast }$ not constant.

Let $M^{3}$ be the manifold $M=\{(x_{1},x_{2},x_{3})\in {\mathbf
R}^{3}\, |\,
x_{3}\neq 0\}$, where $(x_{1},x_{2},x_{3})$ are the standard coordinates on ${\mathbf R}^{3}$. The vector fields
\[
e_{1}=\frac{\partial }{\partial x_{1}},  \qquad
e_{2}=-\ 2x_{2}x_{3}\frac{\partial }{\partial x_{1}}+\frac{2x_{1}}{{x_{3}}^{3}}\frac{\partial }{\partial x_{2}} -\frac{1}{{x_{3}}^{2}}\frac{\partial }{\partial x_{3}}, \qquad
e_{3}=\frac{1}{x_{3}}\frac{\partial }{\partial x_{2}}
\]
are linearly independent at each point of $M$. We get
\[
\lbrack e_{1},e_{2}]=\frac{2}{{x_{3}}^{2}}\,e_{3},\quad \lbrack
e_{2},e_{3}]=2e_{1}+\frac{1}{{x_{3}}^{3}}\,e_{3},\quad \lbrack
e_{3},e_{1}]=0.
\]

Let $g$ be the Riemannian metric defined by $g(e_i,e_j)=\delta_{ij}$, $i,j=1,2,3$, $\nabla$ its Riemannian connection, $R$ the curvature
tensor of $g$ and $\eta $ the 1-form defined by $\eta \left(
X\right) =g(X,e_{1})$, for any $X$ on $M$, which is a contact form
because $\eta \wedge d\eta \neq 0$. Let $\phi $ be the
$(1,1)$-tensor field defined by
\[
\phi e_{1}=0,\quad \phi e_{2}=e_{3},\quad \phi e_{3}=-e_{2}.
\]
Using the linearity of $\phi ,d\eta $ and $g$ we find
\[
\eta (e_{1})=1,\quad \phi ^{2}X=-X+\eta \left( X\right) e_{1},\quad
d\eta \left( X,Y\right) =g(X,\phi Y),
\]
\[
g(\phi X,\phi Y)=g\left( X,Y\right) -\eta \left( X\right) \eta
\left( Y\right),
\]
for any vector fields $X,Y$ on $M$. Hence $(\phi ,e_{1},\eta ,g)$
defines a contact metric structure on $M$. Using Koszul's formula we
obtain:
\[
\begin{array}{ll}
\displaystyle\nabla _{e_{1}}e_{2}=\left(-1+\frac{1}{{x_{3}}^{2}}\right)e_{3}, &\medskip
 \displaystyle\nabla _{e_{2}}e_{1}=-\left(1+\frac{1}{{x_{3}}^{2}}\right)e_{3} \medskip\\
\displaystyle\nabla _{e_{1}}e_{3}=\left( 1-\frac{1}{{x_{3}}^{2}}\right) e_{2}, & \medskip
 \displaystyle\nabla _{e_{3}}e_{1}=\left( 1-\frac{1}{{x_{3}}^{2}}\right) e_{2}, \medskip\\
\displaystyle\nabla _{e_{2}}e_{3}=\left( 1+\frac{1}{{x_{3}}^{2}}\right) e_{1},  & \medskip
 \displaystyle\nabla _{e_{3}}e_{2}=\left( -1+\frac{1}{{x_{3}}^{2}}\right) e_{1}-\frac{1}{{x_3}^3} \ e_{3}, \medskip \\
\displaystyle\nabla _{e_{2}}e_{2}=0, &\displaystyle\nabla_{e_{3}}e_{3}=\frac{1}{{x_{3}}^{3}}\,e_{2}.
\end{array}
\]
The tensor $h$ satisfies
\[
he_{1}=0,\qquad he_{2}=\lambda e_{2},\qquad he_{3}=-\lambda e_{3},
\]
where $\lambda =\displaystyle\frac{1}{{x_{3}}^{2}}$. Putting $\mu =2\left( 1-\displaystyle\frac{1}{{x_{3}}^{2}}\right) $ and
$\kappa =1-\displaystyle\frac{1}{{x_{3}}^{4}}$, we obtain
\[
R\left( X,Y\right) \xi =\kappa \left( \eta \left( Y\right) X-\eta
\left( X\right) Y\right) +\mu \left( \eta \left( Y\right) hX-\eta
\left( X\right) hY\right) .
\]
Therefore $M$ is a generalized $\left( \kappa ,\mu \right) $-space with $\kappa ,\mu $ non-constant functions on $M$.\

Let us now see that the manifold $M$ is also a generalized $\left( \kappa ,\mu \right)$-space form $M(f_{1}^{\ast },0,f_{3}^{\ast },f_{4}^{\ast },0,0)$.
From the definition of Riemannian curvature, we get that the
non-trivial curvatures
are:
\begin{align*}
R(e_{1},e_{2})e_{1} =&-(\kappa +\lambda \mu )e_{2}, \\
R(e_{1},e_{2})e_{2} =&(\kappa + \lambda \mu)e_{1}, \\
R(e_{1},e_{2})e_{3} =&0, \\
R(e_{1},e_{3})e_{1} =&(-\kappa +\lambda \mu )e_{3}, \\
R(e_{1},e_{3})e_{2} =&0, \\
R(e_{1},e_{3})e_{3} =&(\kappa -\lambda \mu )e_{1}, \\
R(e_{2},e_{3})e_{1} =&0, \\
R(e_{2},e_{3})e_{2} =&(\kappa +\mu -2\lambda ^{3})e_{3}, \\
R(e_{2},e_{3})e_{3} =&-(\kappa +\mu -2\lambda ^{3})e_{2}.
\end{align*}

On the other hand, we know that every contact metric generalized
$\left( \kappa ,\mu \right) $-space form $M(f_{1}^{\ast
},0,f_{3}^{\ast
},f_{4}^{\ast },0,0)$ is a generalized $(\kappa ^{\ast },\mu ^{\ast })$-space with $\kappa ^{\ast }=f_{1}^{\ast }-f_{3}^{\ast }$ and $\mu
^{\ast }=f_{4}^{\ast }$. Using (\ref{eq-R-cmgkmsf-3d-1}), we can
write the above curvature values in terms of $\kappa ^{\ast },\mu^{\ast }$ and $f_{1}^{\ast }$. If we make equal both sets of
equations, we obtain a system which can be simplified to:
\[
\left\{
\begin{array}{rcl}
\kappa  & = & \kappa ^{\ast }=f_{1}^{\ast }-f_{3}^{\ast }   \\
\mu  & = & \mu ^{\ast }=f_{4}^{\ast }  \\
f_{1}^{\ast } & = & -\ (\kappa +\mu -2\lambda ^{3}).
\end{array}
\right.
\]

Therefore, the solution to the system is:
\[
\left\{
\begin{array}{rcl}
f_{1}^{\ast } & = & \displaystyle-\ 3+\frac{2}{{x_{3}}^{2}}+\frac{1}{{x_{3}}
^{4}}+\frac{2}{{x_{3}}^{6}} \medskip  \\
f_{3}^{\ast } & = & \displaystyle-\ 4+\frac{2}{{x_{3}}^{2}}+\frac{2}{{x_{3}}
^{4}}+\frac{2}{{x_{3}}^{6}}  \medskip \\
f_{4}^{\ast } & = & \displaystyle 2\left( 1-\frac{1}{{x_{3}}^{2}}\right) .
\end{array}
\right.
\]

Thus we conclude that this example is a contact metric generalized
$\left( \kappa ,\mu \right) $-space form $M^{3}(f_{1}^{\ast
},0,f_{3}^{\ast },f_{4}^{\ast },0,0)$, where $f_{1}^{\ast
},f_{3}^{\ast },f_{4}^{\ast }$ are
non-constant functions. Moreover, its scalar curvature is, by \eqref{eq-tau},
\[
\tau =2(3f_{1}^{\ast }-2f_{3}^{\ast })=2 \left( -1 +\frac{2}{{x_{3}}^{2}} -\frac{1}{{x_{3}}^{4}}+\frac{2}{{x_{3}}^{6}}\right),
\]%
and therefore not constant.

\begin{rem-new}
The second example of generalized $\left( \kappa ,\mu \right)
$-space of dimension $3$ given by T. Kougfogiorgos and C. Tsichlias
in \cite{Kouf-00} is also a contact metric generalized $\left(
\kappa ,\mu \right) $-space form $M^{3}(f_{1}^{\ast},0,f_{3}^{\ast},f_{4}^{\ast},0,0)$ with
non-constant functions:
\[
\left\{
\begin{array}{rcl}
f_{1}^{\ast } & = & \displaystyle-\
3-\frac{2}{{x_{3}}^{4}}+\frac{1}{{x_{3}}^{8}}+\frac{10}{{x_{3}}^{14}} \medskip  \\
f_{3}^{\ast } & = & \displaystyle-\
4-\frac{2}{{x_{3}}^{4}}+\frac{2}{{x_{3}}^{8}}+\frac{10}{{x_{3}}^{14}}  \medskip \\
f_{4}^{\ast } & = & \displaystyle2\left( 1+\frac{1}{{x_{3}}^{4}}\right).
\end{array}
\right.
\]
\end{rem-new}

\section[Trans-Sasakian generalized (k,m)-space forms]{Trans-Sasakian generalized $(\protect\kappa ,\protect\mu )$-space
forms\label{sect-trans}}

We will see in this section that if a manifold is trans-Sasakian, then $h=0$. Hence every generalized $\left( \kappa ,\mu \right) $-space form $M(f_{1},\ldots ,f_{6})$ with a trans-Sasakian structure is a
generalized Sasakian space form (see \cite{ABC,AC-DGA}).

We recall that an almost contact metric manifold $(M,\phi ,\xi ,\eta,g)$
is said to be {\em trans-Sasakian} if there exist functions $\alpha $ and $\beta $ on $M$ such that
\begin{equation}
(\nabla _{X}\phi )Y=\alpha (g\left( X,Y\right) \xi -\eta \left(
Y\right) X)+\beta (g(\phi X,Y)\xi -\eta \left( Y\right) \phi X),
\label{eq-trans-Sas-1}
\end{equation}
for any vector fields $X,Y$ on $M$. In a trans-Sasakian manifold it
is known that
\begin{equation}
\nabla _{X}\xi = -\ \alpha \phi X + \beta \left(X-\eta \left(
X\right) \xi \right).  \label{eq-trans-Sas-2}
\end{equation}

We will now prove two properties that trans-Sasakian manifolds have
in common with contact metric manifolds.

\begin{prop}
\label{prop-trans-Sas-1} If $(M,\phi ,\xi ,\eta ,g)$ is a
trans-Sasakian manifold, then $\nabla _{\xi }\phi =0$ and $\nabla
_{\xi }\xi =0$.
\end{prop}

\begin{proof} It follows directly from (\ref{eq-trans-Sas-1}) and (\ref{eq-trans-Sas-2}). \end{proof}
\begin{prop}
\label{prop-trans-Sas-h} If $(M,\phi,\xi,\eta,g)$ is a
trans-Sasakian manifold, then $h=0$.
\end{prop}

\begin{proof} Using the definition of $h$ and applying the
usual properties we obtain
\[
2hX=(L_\xi \phi)X=(\nabla _{\xi }\phi )X-\nabla _{\phi X}\xi +\phi
\nabla _{X}\xi .
\]
Therefore, the result follows directly from (\ref{eq-trans-Sas-2})
and Proposition~\ref{prop-trans-Sas-1}. \end{proof}
\begin{cor}
Every trans-Sasakian manifold with a contact metric structure is
Sasakian.
\end{cor}

\begin{proof} If $M$ is a trans-Sasakian manifold, it follows that $h=0$ from Proposition~\ref{prop-trans-Sas-h}. If $M$ is also a
contact metric manifold, then it is a $K$-contact manifold and
$\nabla _{X}\xi
=-\phi X$ is satisfied. Comparing such an equation with \eqref{eq-trans-Sas-2}, we deduce that $\alpha =1$ and $\beta =0$.
Substituting these values in (\ref{eq-trans-Sas-1}), we obtain
(\ref{eq-Sas}), which is one of the characterizations of Sasakian
manifolds. \end{proof}

In view of Proposition~\ref{prop-trans-Sas-h}, we have the following

\begin{thm}
Every trans-Sasakian generalized $\left( \kappa ,\mu \right) $-space
form is a generalized Sasakian space form.
\end{thm}

\noindent \textbf{Acknowledgements.} The authors would like to thank the referee for the valuable suggestions.

\end{document}